\magnification\magstep1
\baselineskip = 18pt
\overfullrule = 0pt

\def\n{\noindent}
\def\qed{{\hfill{\vrule height7pt width7pt
depth0pt}\par\bigskip}} 

\def\pf{\medskip\n {\bf Proof.}~~}

\def\EE{{\rm I}\!{\rm E}}
\def\FF{{\rm I}\!{\rm F}}
\def\CC{ \;Ê{}^{ {}_\vert }\!\!\!{\rm C}}

\def\RR{{\rm I}\!{\rm R}}
\def\NN{{\rm I}\!{\rm N}}
\def\nat{\NN}
\def\T{{\bf T}}
\def\ZZ{{\rm Z}\!\!{\rm Z}}
\def\ie{{\it i.e.}}
\def\eg{{\it e.g.}}
\def\tilde{\widetilde}

\def\cf{{\it cf. \ }}
\def\ve{\vfill\eject}

\def\RR{{\mathop{{\rm I}\kern-.2em{\rm R}}\nolimits}}
\def\NN{{\mathop{{\rm I}\kern-.2em{\rm N}}\nolimits}}
\def\FF{{\mathop{{\rm I}\kern-.2em{\rm F}}\nolimits}}
\def\EE{{\mathop{{\rm I}\kern-.2em{\rm E}}\nolimits}}
\def\CC{{\rm C\kern-.18cm\vrule width.6pt height 6pt depth-.2pt
\kern.18cm}}
\def\ZZ{{\mathop{{\rm Z}\kern-.28em{\rm Z}}\nolimits}}

\def\nat{\NN}

\def\PP{{\mathop{{\rm I}\kern-.2em{\rm P}}\nolimits}}
\def\RR{{\mathop{{\rm I}\kern-.2em{\rm R}}\nolimits}}

\centerline{\bf Multipliers of the Hardy space $H^1$ and}
\centerline{\bf power bounded operators}

\centerline{by }
\centerline{ Gilles Pisier\footnote*{{
Supported in part by NSF   and
by the Texas Advanced Research Program 010366-163.}}}
\centerline{ Texas A\&M University}
\centerline{College Station, TX 77843, U. S. A.}
\centerline{and}
\centerline{Universit\'e Paris VI}
\centerline{Equipe d'Analyse, Case 186, 75252}
\centerline{ Paris Cedex 05, France}\bigskip

\noindent {\bf Abstract}
We study the space of functions $\varphi\colon \ \NN\to  \CC$ such that there
is a Hilbert space $H$, a power
 bounded operator $T$ in $B(H)$
 and vectors
$\xi,\eta$ in $H$ such that
$$\varphi(n) = \langle T^n\xi,\eta\rangle.$$
This implies that the matrix $(\varphi(i+j))_{i,j\ge  0}$ is a Schur
multiplier of $B(\ell_2)$ or equivalently is in the space $(\ell_1
\buildrel {\vee}\over {\otimes}  \ell_1)^*$. We show
that the converse does not hold, which answers a question
raised by Peller [Pe]. Our approach makes use of a new class
of Fourier multipliers of $H^1$ which we call
``shift-bounded''. We show that there is a $\varphi$ which is
a ``completely bounded'' multiplier of $H^1$, or equivalently
for which
$(\varphi(i+j))_{i,j\ge 0}$ is a bounded Schur multiplier of $B(\ell_2)$, but
which is not ``shift-bounded'' on $H^1$. We also give a characterization of
``completely shift-bounded'' multipliers on $H^1$.
\vskip.3in


\ve
\n {\bf \S 0. Introduction} 

\n The main motivation of this paper
is a question of Peller on power bounded operators [Pe]. To
state it, we need some specific notation, as follows. 

\n For any
$c>1$ and any polynomial $P(z) = \sum a_nz^n$, let 
$$  |||P|||_c = \sup \{ \| \sum a_n T^n\| \}$$
where the
supremum runs over all (power bounded) operators $T$ in
$B(\ell_2)$ such that $\sup\limits_{n\ge 1} \|T^n\| \le c$. We also let
$$\|P\|_{\cal L} = \inf\{\|A\|_{\ell_1\buildrel {\vee}\over
{\otimes}\ell_1}\}$$ where the infimum runs over all elements
$A = \sum A_{ij} e_i\otimes e_j$ in the injective tensor product
 $\ell_1\buildrel {\vee}\over {\otimes}\ell_1$ such that
$$a_n = \sum_{i+j=n} A_{ij}.$$
As observed by Peller, it follows from Grothendieck's theorem that there is a
constant $K$ such that for all $c>1$ we have
$$|||P|||_c \le Kc^2\|P\|_{\cal L}.$$
Peller asked whether conversely there is any $c>1$ for which
$|||~~|||_c$ and
$\|~~\|_{\cal L}$ are equivalent. We will prove below that it is not so.
Unfortunately, Peller's basic question whether all the norms $|||~~|||_c$ are
equivalent for all $c>1$ remains open, although we propose a ``new'' approach
for its solution, directly inspired by Peller's ideas
 in [Pe] but
revised in  light of the recently developed operator
 space theory (see especially [BRS, B, BP2]). As we show
below, the latter theory clearly suggests that one 
should replace the iterated
injective tensor products $  \ell_1\buildrel {\vee}\over
{\otimes}  \cdots  \buildrel {\vee}\over
{\otimes}\ell_1  $ ($d$ times), which Peller uses, by the
iterated Haagerup tensor products $\ell_1
\otimes_h \cdots \otimes_h \ell_1$. (Actually, since $\NN$ is 
commutative, we should consider the
``symmetrized" iterated Haagerup tensor products of [OP], but
this can be left implicit in this note.)

Peller observed that for any two polynomials $P,Q$ we have
$$|||PQ|||_c \le |||P|||_c \ |||Q|||_c.$$
In other words, $|||~~|||_c$ is a Banach algebra norm. Moreover, its
definition clearly shows that it is an ``operator algebra'' norm, i.e.\ the
resulting Banach algebra can be isometrically embedded into the algebra $B(H)$
of all bounded operators on a Hilbert space $H$. Thus Peller was
led to ask whether $\|~~\|_{\cal L}$ is equivalent to an
operator algebra norm, or merely even to a Banach algebra one. 
In his review of Peller's paper 
  ([Math. Reviews Sept.
1983i,  47019]),
G.~Bennett proved that $\|~~\|_{\cal L}$ is indeed a Banach algebra norm, but
we will prove below 
(see Corollary 2.2) that it is not equivalent to an operator
algebra norm. As Peller observed, it suffices to prove that for
any
$c>1$  the norms
$\|~~\|_{\cal L}$ and $|||~~|||_c$  are not equivalent. 

\n To prove this, we use a connection with bounded (Fourier)
multipliers on the Hardy space $H^1$, as follows. Given a
Banach space $B$, we denote by $H^1(B)$ the completion the
$B$-valued polynomials $f(z) = \sum x_nz^n$ (here $x_n\in B$)
for the norm
$$\|f\|_{H^1(B)} = \int \|f(z)\|_B \ dm(z)$$
where $m$ denotes the normalized Lebesgue measure of the unit circle. 

\n When $B =\CC$, this is the classical Hardy space, which we
simply denote as usual by $H^1$.

 Let $\varphi\colon \ \NN\to
\CC$ be a function in ${\cal L}^*$, i.e.\ we assume there is a
constant $C$ such that for any polynomial
$P(z) = \sum a_nz^n$ (with $a_n\in \CC$) we have
$$\left|\sum \varphi(n)a_n\right| \le C\left\|\sum a_nz^n\right\|_{\cal
L}.\leqno (0.1)$$
Then $\varphi$ defines a multiplier
$$M_\varphi\colon \ \sum a_nz^n \to \sum a_n \varphi(n)z^n$$
which is bounded on $H^1$.
Actually, $M_\varphi$ is ``completely bounded'' on $H^1$
 (see [P1, \S 6]) which
means that $M_\varphi$ defines a bounded multiplier on $H^1(S_1)$ where $S_1$
denotes the Banach space of all trace class operators on $\ell_2$ (equipped
with the norm $\|x\|_{S_1} = {\rm tr}(|x|)$). Conversely, as
observed in [P1, Th. 6.2] any
$\varphi$ such that $M_\varphi$ is completely bounded on $H^1$ is in ${\cal
L}^*$, \ie\ satisfies (0.1), and the norms $\|\varphi\|_{{\cal
L}^*}$ and
$\|M_\varphi\colon  \ H^1(S_1)\to H^1(S_1)\|$ are equivalent
(see (1.5) and (1.6) below). The main virtue of this note is the
introduction of a ``restricted'' class of completely bounded
(in short c.b.) multipliers on
$H^1$. We will say 
that $M_\varphi$ is ``shift-bounded'' on $H^1$ if for any $x =
\sum\limits_{n\ge 0} x_nz^n$ in $H^1$ we have
$$\int \sup_{k\ge 0} \left|M_\varphi(z^k x)\right| dm(z)=
\int \sup_{k\ge 0} \left|\sum_{n\ge
0} x_n\varphi(n+k)z^n\right| dm(z) <
\infty.$$
We will show below that any $\varphi$ such that $|||\varphi|||^*_c < \infty$
for some $c>1$ must define a ``shift-bounded'' multiplier on
$H^1$. Thus to show that $\|~~\|_{\cal L}$ and $|||~~|||_c$ are
not equivalent it suffices to produce a multiplier on $H^1$
which is completely bounded but {\it not\/} shift-bounded.
This follows from our main result (see Theorem 2.1 below).
In \S 3 we include some remarks on the class of shift-bounded
multipliers on $H^1$ and a characterization of their  
analogue on
$S_1$-valued
$H^1$, the "completely shift-bounded'' ones, 
 which might  be of independent interest. Although this uses
ideas and techniques from the
recently developed "operator space theory'' ([BP2,ER1-3]),
our formulation  (especially in Theorem 3.3)
hopefully will be accessible to readers not familiar with it.

\n {\bf \S 1. Notation and  background}

Let $G$ be a 
semi-group with unit. Since we mostly concentrate here on the
case $G=\nat$, we will denote the operation
of $G$ additively with unit $0$. However, most of
our notation makes sense for a general (non-commutative)
semi-group. In particular we refer to [P5] for
a detailed treatment of the case when $G$ is a free group
(see also [P2] for related results).
After some hesitation, we chose to write separate papers since,
although multipliers appear in both of these
two cases, 
the same questions   require quite different techniques.

Let $\pi\colon
\ G\to B(H)$  be a uniformly bounded   
unital semi-group homomorphism, \ie\  we have
$$\pi(s+t)= \pi(s)\pi(t) \quad \pi(0)=I.$$
We denote
$$|\pi| = \sup \{\|\pi(t)\|_{B(H)}\mid t\in G\}.$$
 
 Let $c\ge 1$ .  We denote by
$B_c(G)$ the space of all ``matrix coefficients" 
of the   
unital semi-group homomorphisms which are
uniformly bounded by $c$. More precisely,  $B_c(G)$ is the space
of functions 
${\varphi}\colon \ G\to \CC$ for which 
there is    
$\pi\colon
\ G\to B(H)$ as above with 
$|\pi|\le c$ together with vectors $\xi,\eta$ in $H$ such that
$${\varphi}(t) = \langle \pi(t)\xi,\eta\rangle.\leqno
(1.1)\ \forall~t\in G$$   Moreover, we denote
$$\|{\varphi}\|_{B_c(G)} = \inf\{\|\xi\|\ \|\eta\|\mid {\varphi}(\cdot) = \langle\pi(\cdot)\xi, 
\eta\rangle \hbox{ with } |\pi|\le c\}.$$
Note that when $c=1$ and $G$ is a group,
 $B_1(G)$ coincides with the classical space of
coefficients of unitary representations of $G$, usually
denoted  by
   $B(G)$, with the same norm.
  Indeed, it is easy to 
check in the group case that
$|\pi|=1$ iff $\pi$ is a unitary  
representation.

The space $B_c(G)$ is a Banach space (for the above  norm). Moreover, for any 
$c'\ge 1$ we have
$$f\in B_c(G), g\in B_{c'}(G)\Rightarrow f\cdot g \in B_{cc'}(G).$$
Note moreover that if $c\ge c'$ we have a norm one inclusion
$$B_{c'}(G) \subset B_c(G).$$
 
In the main case of interest to us here, $G=\NN$, we have
an isometric identity
$$B_1(\NN)= A(D)^*\leqno(1.2)$$
   Here $A(D)$ is the disc algebra which can be
defined as the completion of the space of (analytic) polynomials
$P$ under the sup-norm over the unit disc in $\CC$. Indeed, by
a well known inequality
of von Neumann (see \eg \ [P1, \S 1]) for any such $P$
and for any contraction $T$ in $B(H)$ (meaning $\|T\|\le 1$), 
we have 
$$\|P(T)\| \le \|P\|_{A(D)}.$$
A unital homomorphism $\pi\colon
\ \NN\to B(H)$ is in $1-1$ correspondence with a   $T\in B(H)$
such that
$\forall n\in \NN\quad\pi(n)=T^n,$ thus,
 for any ${\varphi}\colon \ \NN\to \CC$, we have
$\|{\varphi}\|_{B_1(\NN)}\le 1$ iff there is a contraction
$T$ and $\xi,\eta$ in the unit ball of $H$ such that
$${\varphi}(n)=\langle T^n \xi, \eta \rangle,\leqno(1.3)$$
and by von Neumann's inequality this holds
iff  $\|{\varphi}\|_{ A(D)^*}\le 1$.
This verifies (1.2).
More generallly, for any $c\ge 1$, we have 
$\|{\varphi}\|_{B_c(\NN)}\le 1$ iff
there is a power bounded $T\in B(H)$ with
$\sup_n \|T^n\|\le c$ and $\xi,\eta$ in the unit ball of $H$
such that  (1.3) holds.

Let $d\ge 1$ be an integer.
  Let
$M_d(G)$ be
 the space of all 
functions ${\varphi}\colon \ G\to \CC$ such that there are bounded functions 
$\xi_i\colon \ G\to B(H_i, H_{i-1})$ ($H_i$ Hilbert) with $H_0 = \CC$, $H_d 
=\CC$ such that
$${\varphi}(t_1+t_2+\ldots +t_d) = 
\xi_1(t_1) \xi_2(t_2)\ldots \xi_d(t_d).\leqno (1.4)\qquad 
\forall~t_i\in G$$
Here of course we use the identification $B(H_0,H_d) = B(\CC,\CC)\simeq \CC$. We 
define
$$\|{\varphi}\|_{M_d(G)} = \inf\{\sup_{t_1\in G} \|\xi_1(t_1)\|\ldots \sup_{t_d\in G} 
\|\xi_d(t_d)\|\}$$
where the infimum runs over all
 possible ways to write ${\varphi}$ as in (1.4).

\n  It is quite easy to see that $M_d(G)$
is a Banach algebra for the pointwise product of
functions on $G$. 

When $G=\NN$,   we have
$\|\varphi\|_{M_2(\NN)}\le 1$ iff there are  sequences $(x_n)$
and $(y_n)$ in the unit ball of $H$
such that
$$\forall i,j\in \NN \qquad \varphi(i+j)=\langle x_i, y_j\rangle.$$
Equivalently, if we denote by
$u_\varphi \colon \ell_1 \to \ell _\infty$ the linear
operator
with (Hankel) matrix \break $(\varphi(i+j))$, we have 
$\|\varphi\|_{M_2(\NN)}=\gamma_2(u_\varphi)$ (here 
$\gamma_2(.)$ is the norm of factorization through a
Hilbert space). 
A fortiori we have
$$K^{-1}\|\varphi\|_{{\cal L}^*} \le \|\varphi\|_{M_2(\NN)} \le 
\|\varphi\|_{{\cal L}^*}\leqno(1.5)$$ where $K$ denotes the
Grothendieck constant.
Moreover, we have
$$ \|\varphi\|_{M_2(\NN)}=\|M_\varphi\colon  \ H^1(S_1)\to
H^1(S_1)\|\ge \|M_\varphi\colon  \ H^1 \to
H^1 \|.\leqno(1.6)$$ 
 For example, any bounded $\varphi$ with support in a lacunary
sequence such as 
$\{ 2^n\mid n\ge 0\}$ is in $M_2(\NN)$, see Lemma 2.4  below for
a more general fact.

\n   See
[Bo]  and [P1,
\S 5-6]  for more on all this.

The definition of the spaces ${M_d(G)}$ (and of
the  $B(H)$-valued version of these spaces for which we refer
to [P5]) is motivated by the work of Christensen-Sinclair on
''completely bounded multilinear maps" and the so-called
Haagerup tensor product (see [CS]). The
connection is explained in detail in [P3, P5],
and is important for the results below, but we prefer
to skip this in the present, hopefully more accessible, 
exposition.

  Note the following easily
checked
  inclusions,  valid when $G$ is any
semigroup with unit:
$$\eqalign{B(G)=B_1(G) &\subset \bigcup_{c>1}
B_c(G)
\subset M_d(G) \subset  M_{d-1}(G) 
\subset\cdots\cr
\cdots &\subset M_2(G) \subset M_1(G) =
\ell_\infty(G),}$$ and we have
$$\|f\|_{M_{m}(G)} \le
\|f\|_{M_d(G)}.\leqno
(1.7)\qquad \forall m\le d$$
Moreover, we have
$$\|{\varphi}\|_{M_d(G)} \le c^d\|{\varphi}\|_{B_c(G)}.\leqno
(1.8)\qquad \forall~{\varphi}\in B_c(G)$$ Indeed, if
${\varphi}(\cdot) = \langle\pi(\cdot)\xi,\eta\rangle$ with
$|\pi|\le c$, then  we can write
$$\eqalign{{\varphi}(t_1+t_2+\ldots +t_d) &= \langle \pi(t_1)\ldots 
\pi(t_d)\xi,\eta\rangle\cr
&= \xi_1(t_1)\xi_2(t_2)\ldots \xi_d(t_d)}$$
where $\xi_1(t_1) \in B(H_\pi,\CC)$, $\xi_d(t_d) \in B(\CC,H_\pi)$ and 
$\xi_i(t_i) \in B(H_\pi,H_\pi)$ $(1<i<d)$ are defined by $\xi_1(t_1)h = 
\langle\pi_1(t_1)h,\eta\rangle$ $(h\in H_\pi)$ $\xi_d(t_d)\lambda = 
\lambda\pi(t_d)\xi$ $(\lambda\in\CC)$ and $\xi_i(t_i)= \pi(t_i)$ $(1<i<d)$. 
Therefore, we have
$$\eqalign{\|{\varphi}\|_{M_d(G)} &\le \sup\|\xi_1\| \sup\|\xi_2\|\ldots\sup 
\|\xi_d\|\cr
&\le |\pi|^d\|\xi\|\ \|\eta\|\le c^d\|\xi\|\ \|\eta\|}$$
whence the announced inequality (1.8).

 Let $UB(G) =
\bigcup\limits_{c>1} B_c(G)$.
Then ${\varphi}\in UB(G)$ iff 
$\sup\limits_{m\ge 1} \|{\varphi}\|^{1/m}_{M_m(G)} <
\infty$. More precisely, let $c({\varphi})$ denote the
infimum of the numbers $c\ge 1$ for which
${\varphi}\in B_c(G)$. Then, we have (see [P5])
$$c({\varphi})=\limsup\limits_{m\to \infty}
\|{\varphi}\|^{1/m}_{M_m(G)}.$$
Moreover, it follows from [BP2] (see also [P4, \S 7])
that $$\|{\varphi}\|_{B_1(G)} =\sup\limits_{m}
\|{\varphi}\|_{M_m(G)}.$$

  The definition of the spaces $B_c(G)$ and
$M_d(G)$ shows that they are dual spaces.
There is a natural duality between these spaces
and the semi-group algebra $\CC[G]$ which we view as
the convolution algebra of finitely supported functions on $G$.
Indeed, for any function $f\colon\ G\to \CC$
and any  $g$ in $\CC[G]$, we set
$$<g,f>= \sum_{t\in G} g(t)f(t),$$
and we define the spaces $X_d(G)$ and $\tilde A_c$
respectively as the completion of $\CC[G]$
for the respective norms
$$\|g\|_{X_d(G)}= \sup \{ |<g,f>| \mid f\in M_d(G), \
\|f\|_{M_d(G)}\le 1
\} $$
and
$$ |||g|||_c= \sup \{ |<g,f>| \mid f\in B_c(G), \
\|f\|_{B_c(G)}\le 1
\} .$$
 Then the following isometric identities   are rather easy to
check: 
$$ {B_c(G)}= (\tilde A_c)^*\quad{\rm and}\quad
M_d(G)=(X_d(G))^*.$$ Obviously, we can also write (here we can
restrict to
 $H=\ell_2$ if we wish)
$$ |||g|||_c = \sup \{ \|\sum g(t) \pi(t)\|\mid
\pi\colon \ G\to B(H),\  |\pi|\le c
\} .$$
Thus, when $G=\NN$, ${\tilde A_c}$  can be identified with the
completion of polynomials (here a polynomial is identified
with the sequence of its coefficients) 
for the norm $ |||~~|||_c~$ introduced in \S 0. 
  The last formula shows that ${\tilde
A_c}$ is naturally equipped with
an operator algebra structure under convolution: we have
$|||g_1 *g_2|||_c\le |||g_1  |||_c ||| 
 g_2|||_c.$

\n
However, the analogue for the spaces $X_d(G)$ fails in general.
(This is   the basic idea
used by Haagerup   to prove that
$M_2(\FF_\infty)\not=B_c(\FF_\infty)$ for any $c$
(see Remark 1.2 below)\ :
he proves first in [H] that
$X_2(\FF_\infty)$ is not a Banach algebra under convolution.)

In sharp contrast, as we already mentioned, $X_2(\NN)$ is indeed
a Banach algebra!

\n Although $X_d(G)$  is not in general
a Banach algebra  under convolution, it satisfies the following
property: if $g_1\in X_d(G)$ and $g_2\in X_k(G)$,
then $g_1 *g_2\in X_{d+k}(G)$ and
$$\|g_1 *g_2\|_{ X_{d+k}(G)}\le \|g_1  \|_{X_d(G)} \| 
 g_2\|_{X_k(G)}.\leqno(1.9)$$
 See [P5] for a detailed  proof.

To explain the relevance of the spaces $M_d(G)$ 
for Peller's question, we quote

\proclaim Theorem 1.1. ([P4, P5]). Let $G$ be a
semigroup with unit. 
The following assertions are
equivalent:
\item {\rm (i)} There is a $\theta\ge 1$ such
that 
$B_{\theta}(G) =
B_c(G)$ for all $c>\theta$.
\item {\rm (i)'} There is a $\theta\ge 1$ such
that 
$B_{\theta}(G) =
B_c(G)$ for some $c>\theta$.
\item {\rm (ii)}
 There are   $\theta\ge 1$ and  
an integer $d$ such that
$B_{\theta}(G) =M_d(G)$.
\item {\rm (iii)} There is
an integer $d$ such that  $ M_d(G)=M_{2d}(G)$.
\item {\rm (iv)} There is
an integer $d$ such that $ X_d(G)$ is (up to isomorphism)
a unital operator algebra under convolution.

Thus to  show that the   norms
$ |||~~|||_c$ and $ |||~~|||_\theta$ are  not equivalent
whenever  
$c\not=\theta$, it suffices to prove the following

\proclaim Conjecture. 
$ M_d(\NN)\not=
M_{d+1}(\NN) \quad \forall d\ge 1.$

\n {\bf Remark 1.2} When $d=2$, and $G$ is a group, the space
$M_2(G)$ is the classical space of ``Herz-Schur 
multipliers'' on
$G$. This space also coincides (see  [BoF] or
[P2, p. 110]) with the  space of all c.b.\
``Fourier multipliers'' on the reduced
$C^*$-algebra 
$C^*_\lambda(G)$. The question whether 
 $M_2(G)=UB(G)$   remained open for a while but
Haagerup [H] showed that it is not the  case. More precisely,
he showed that if $G = \FF_\infty$, we have
$$B_c(G) \mathop{\subset}\limits_{\ne}
M_2(G).\leqno\forall~c>1$$
In [P5] we give a different proof of this. 
 More generally we show there that
if $G$ is a non-commutative  free
group, for any $d\ge 1$, we have
$$M_d(G)\not= M_{d+1}(G),$$ and hence there are
elements of
$M_d(G)$  which are  not coefficients
of uniformly bounded representations.

 \bigskip

\n {\bf \S 2. Main results}

\n We wish to prove here a special case of this conjecture,
which answers a question of Peller in [Pe].

\proclaim Theorem 2.1.
$$M_2(\NN)\not=M_3(\NN).$$

%\proclaim Corollary. ${\cal L}=X_2(\nat)$  is a Banach
%algebra (by G. Bennett's aforementioned result)
%but   not an operator algebra 
%under convolution. Equivalently, there is no
%$c$ for which the norms $ |||~~|||_c$ and $\|~~\|_{\cal L}$
%are equivalent. 

\proclaim Corollary 2.2.  For any $c>1$ the norms
$\|~~\|_{\cal L}$ and $|||~~|||_c$  are not equivalent. 
More generally, $\|~~\|_{\cal L}$ is not equivalent to any 
  operator algebra norm.

\pf Note that by (1.5) we have for all polynomial $P$
$$ K^{-1} \| P\|_{X_2(\NN)} \le \| P\|_{\cal L} \le \|
P\|_{X_2(\NN)} .$$ 
So if $\|~~\|_{\cal L}$ and $|||~~|||_c$  were equivalent,
we would have (by duality) $M_2 =B_c $, hence a fortiori
 $M_2=M_3$ which contradicts Theorem 2.1. The second
assertion follows from Theorem 1.1.\qed

\n {\bf Remark.}   This corollary is closely related
to the more recent result due to Kalton and Le
Merdy [KLM] asserting that, for any
$c>1$, there are  power bounded operators which
are not similar to operators with powers bounded
by $c$. Indeed, 
their result implies   the operator space
(=completely bounded) analogue
 of the inequivalence of $ |||~~|||_c$ and  $\|~~\|_{X_d(\nat)}$
(hence in particular of $ |||~~|||_c$ and $\|~~\|_{X_2(\nat)}$)
for all $c>1$ and $d\ge 1$. It also implies the
 completely bounded analogue of the above conjecture.
Equivalently, this shows that for any $d$ there is a
$B(H)$-valued function which is in 
the   operator valued analogue of $M_d(\nat)$
but not in  the  corresponding analogue of $M_{d+1}(\nat)$.

We will consider the Hardy spaces $H^p$. We define $H^p$ as the subspace of all 
functions $x\in L_p(T,m)$ such that the Fourier transform $\hat x\colon \ \ZZ\to 
\CC$ vanishes on the negative integers. We write 
abusively $x = \sum  
\hat x(n)z^n$, meaning that $x$ admits $\sum \hat x(n)z^n$ as its  formal 
Fourier series.

\proclaim Lemma 2.3. Any ${\varphi}$ in $M_3(\NN)$ defines a
shift-bounded multiplier on $H^1$.  More precisely,
 we
have
$$\sup_{x\in B_{H^1}}\left\{\int \sup_{k\ge 0}
 \left|\sum_{n\ge 0} \hat
x(n)z^n {\varphi}(n+k)\right| dm \ 
\right\} \le \|{\varphi}||_{M_3(\NN)}.\leqno (2.1)$$

\pf Note that for each fixed integer $k$ the function $n\to {\varphi}(n+k)$ is in 
$M_2(\NN)$ hence defines a bounded multiplier on $H^1$. Therefore, for any $x$ 
in $H^1$, the series $\sum\limits_{n\ge 0} \hat x(n)z^n {\varphi}(n+~k)$ is in $H^1$, so 
    (2.1) expresses a sort of uniform boundedness
of this family of multipliers. Now assume
$\|{\varphi}\|_{M_3(\NN)}<1$ so that there  are $\xi_i
\in {\ell_2}^*$, $\eta_j\in \ell_2$, $T_k\in
B(\ell_2)$ with
$$\max(\|\xi_i \|, \|\eta_j\|, \|T_k\|) < 1\leqno (2.2)$$
such that
$${\varphi}(i+k+j)  = \langle \xi_i, T_k\eta_j\rangle.\leqno
(2.3)$$ Let $x\in B_{H^1}$. Then, by a classical 
result (\cf \eg \ [Ga p.
87] or [Ni]) we can write $x=gh$ with 
$g,h \in B_{H^2}$.
 Then 
$$\hat x(n) = \sum_{i+j=n} \hat g(i) \hat h(j).$$
Hence
$$\eqalign{\sum_n \hat x(n) z^n  {\varphi}(n+k) &= \sum_{i,j} \hat g(i) \hat h(j) 
z^{i+j} \langle \xi_i, T_k\eta_j\rangle\cr
&= \langle G(z), T_kH(z)\rangle}$$
where
$G(z) = \sum\limits_i \hat g(i) z^i \xi_i \in H^2({\ell_2}^*)$
and
$H(z) = \sum\limits_j \hat h(j) z^j \eta_j \in H^2(\ell_2)$
with
$$\|G\|_{H^2({\ell_2}^*)} = \left(\sum |\hat g(i)|^2
\|\xi_i\|^2\right)^{1/2} \le  1$$
and similarly $\|H\|_{H^2(\ell_2)} \le 1$. Thus we obtain (by Cauchy-Schwarz)
$$\eqalign{\int \sup_k \left|\sum_n \hat x(n)z^n  
{\varphi}(n+k)\right| dm &\le 
\sup_k \|T_k\| \int \|G(z)\|_{{\ell_2}^*} \|H(z)\|_{\ell_2}
dm(z)\cr &\le \|G\|_{H^2({\ell_2}^*)} \|H\|_{H^2(\ell_2)} \le
1.}$$ \qed

\n {\bf Remark.} Actually (see Theorem 3.1 below), it is
possible  to show that $\|{\varphi}\|_{M_3}$ coincides 
with the c.b.\ norm of the ``multiplier'' defined by ${\varphi}$ as above but acting 
from $H^1$ to $H^1(\ell_\infty)$. 

We will use the following well known lemma (see \eg \ [Bo]).

\proclaim Lemma 2.4. For any ${\varphi}\colon \ \NN\to \CC$ we have
$$\|{\varphi}\|_{M_2(\NN)} \le 4 \sup_{n\ge 0}
\left(\sum_{2^n\le i < 2^{n+1}} 
|{\varphi}(i)|^2\right)^{1/2} + |{\varphi}(0)| .$$

\pf Let
$(e_i)_{i\ge 0}$ be the canonical basis of $\ell_2$.
We can write
$${\varphi}(i+j) = \langle x_i,e_j\rangle + \langle x_j,e_i\rangle +
\langle y_i,e_j\rangle$$
where  $x_i = \sum\limits_{i\le k<2i} {\varphi}(k) e_{k-i}$ 
and $y_i = {\varphi}(2i)e_i$.  Therefore
$$\|{\varphi}\|_{M_2(\NN)} \le 2\sup_{i\ge 0} \|x_i\| +
\sup_i \|y_i\|.$$ Let $C = \sup\limits_{n\ge 0}
\left(\sum\limits_{2^n\le i<2^{n+1}}
|{\varphi}(i)|^2\right)^{1/2}$. We have $x_0 = 0$,
$\|x_i\| \le\sqrt 2 \ C$ and
$\|y_i\| \le C + |{\varphi}(0)|$. Hence
$$\|{\varphi}\|_{M_2(\NN)} \le (2\sqrt 2 + 1) C +
|{\varphi}(0)|.$$
\qed

\proclaim Lemma 2.5. Let $F_1,\ldots, F_K$ be analytic trigonometric
polynomials with degree $\le 2^{K-1}$. Then let
$\varphi =\hat F$ with 
$$F = \sum_{K/2 < p < K}  z^{2^{2p}} F_p.$$
We have then:
$$\int \sup_{K/2 < p < K} |F_p| dm \le 3 \|\varphi\|_{M_3(\NN)}.$$

\pf Let $x\in H^1$ be the  (La Vall\'ee Poussin type)
 kernel such that: $\hat x\equiv 1$ on the interval
$[2^K,2^{K+1}]$,
$\hat x(0)=0$,  $\hat x\equiv 0$ on the interval $[
3.2^K, \infty($ and $\hat x$ is linear on the
remaining intervals $[0,2^K]$ and $[2^{K+1}, 3.
2^K]$. A well known computation shows
 that $\|x\|_1 \le 2$. By
Lemma~2.3 we have
$$\int \sup_{k\ge 0} \left|\sum_{n\ge k} \hat x(n-k) \varphi(n) z^n\right| dm
\le 2\|\varphi\|_{M_3(\NN)}.$$
For each $p$ with $K < 2p<2K$ we let $k(p) = 2^{2p} -2^K$.

\n Hence we have
$$\int \sup_{K/2<p<K} \left|\sum_{n\ge k(p)} \hat x(n-k(p))
\varphi(n)z^n\right| dm \le 2\|\varphi\|_{M_3(\NN)}.$$
Let $A_p = \{n\mid \hat x(n-k(p))\ne 0\}$. We have $A_p\subset k(p) +
[0, 3.2^K]$ hence $A_p \subset [2^{2p}-2^{K}, 2^{2p}+2.2^K]$.
Therefore (since $\varphi$ is supported in the union
of the intervals $[2^{2p}, 2^{2p}+2^{K-1}]$) we find
$$A_p\cap \{n\mid \varphi(n)\ne 0\}
 \subset [2^{2p}, 2^{2p}+2^{K-1}].$$
Now since
$$[2^{2p}, 2^{2p}+ 2^{K-1}] - k(p) \subset [2^K,
2^{K+1}]$$
 we have
$\hat x(n-k(p)) = 1$ for all $n\in A_p\cap \{n\mid \varphi(n)\ne 0\}$. Since
$A_p\cap \{n\mid \varphi(n)\ne 0\} 
\subset [2^{2p}, 2^{2p}+2^{K-1}]$ we must have
simply
$$\sum_{n\ge k(p)} \hat x(n-k(p)) \varphi(n) z^n = z^{2^{2p}} F_p$$
and we conclude that
$$\int \sup_{K/2 < p < K} |F_p| dm \le 2\|\varphi\|_{M_3(\NN)}.$$
\qed

Let $q$ be an integer and let ${\cal P}(q)$ denote the space of all analytic
trigonometric polynomials with degree at most $q$. We define
$$C(q) =\sup\left\{\int \sup_{1\le p\le q} |F_p| dm \ \mid\  F_p\in {\cal P}(q),
\ \sup_{1\le p\le q} \|F_p\|_2 \le 1\right\}.$$

\proclaim Lemma 2.6. For each even integer $K>1$ there is a function
$\varphi_K\colon \ \NN\to \CC$ with support in $[0,2^{2K}]$ such that
$\|\varphi_K\|_{M_2(\NN)} \le 1$ but such that
$$\|\varphi_K\|_{M_3(\NN)} \ge (1/8) C(K/2-1).$$

\pf Let $q = K/2-1$. Let $F_1,\ldots, F_q\in {\cal P}(q)$ be such that
$\sup\limits_{p\le q} \|F_p\|_2 \le 1$ and
$$\int \sup_{p\le q} |F_p| dm = C(q).$$
We consider the function
$$F = \sum_{K/2 < p < K} z^{2^{2p}}F_{p-K/2},$$
and we let $\varphi = \widehat F$.

\n Then by Lemma 2.5 (note that $q=K/2-1\le 2^{K-1}$) we have
$$C(q) \le 2 \|\varphi\|_{M_3(\NN)},$$
and on the other hand by Lemma 2.4 we have
$$\|\varphi\|_{M_2(\NN)} \le 4 \sup_{p\le q}\|F_p\|_2\le 4,$$
whence $\varphi_K = \varphi/4$ has the announced property.\qed

The following fact is elementary and well known:

\proclaim Lemma 2.7. There is a number $\delta>0$ such that for any $q\ge 1$
$$\delta\sqrt q \le C(q) \le \sqrt q.$$

\pf The upper bound is an easy exercise (using $\sup_p |F_p| \le
(\sum |F_p|^2)^{1/2}$). For the lower bound, we use the following well known
consequence of S.~Bernstein's inequality: \ There is a finite subset $B_q\subset {\T}$
such that for any function $F$ in ${\cal P}(q)$ we have
$$\|F\|_\infty \le a \sup_{t\in B_q} |F(t)|$$
and moreover such that
$$|B_q| \le b q$$
where $a\ge 1$ and $b \ge 1$ are absolute constants. We then set $S(z) =
\sum^q_0 z^i$ and for each $\xi$ in $B_q$ we set
$$F_\xi(z) = S(\xi z).$$
We have then for any $z$ in ${\T}$
$$\sup_{\xi\in B_q} |F_\xi(z)|\ge (q+1)/a$$
hence $\int \sup\limits_{\xi\in B_q} |F_\xi(z)|dm(z) \ge (q+1)/a$, and on the
other hand
$$\sup_{\xi\in B_q} \|F_\xi\|_2  \le (q+1)^{1/2}.$$
Hence (assuming without loss of generality that $b$ is an integer $\ge 1$)
this shows that
$$C(b q) \ge (q+1)^{1/2}/a.$$
On the other hand it is easy to check that $C(bq) \le b C(q)$, hence we
finally obtain the announced result with $\delta = (ab)^{-1}$.\qed

\n {\bf Proof of Theorem 2.1.} From the preceding two lemmas, it is clear that 
the norms of $M_2(\NN)$ and $M_3(\NN)$ are not equivalent,
whence Theorem~2.1.\qed

To reformulate more precisely Theorem 2.1, we
need more notation: \ For any integer $d\ge 2$ we denote for
any $\varphi\colon \ \NN\to
\CC$
$$\|\varphi\|_{[d]} = \sup\{\|\varphi\psi\|_{M_d(\NN)}\}$$
where the supremum runs over all $\psi$ in $\ell_\infty(\NN)$ with
$\|\psi\|_\infty\le 1$.
It is well known that  $\|\varphi\|_{[2]}<\infty$
iff $M_\varphi$
maps $H^1$ boundedly into $H^2$, whence (see [Ru, \S 8.6], see
also [Bo]) the following characterization:
$$4^{-1}\|\varphi\|_{[2]}\le  |\varphi(0)| + \sup_{n\ge
0}
\left(\sum_{2^n\le k<2^{n+1}} |\varphi(k)|^2
\right)^{1/2}\le 4 \|\varphi\|_{[2]}.$$
(The left side follows from Lemma 2.4, the other one from
a routine averaging argument and Khintchine's inequality  which
  show that if $\|\varphi\|_{[2]}<\infty$ then
$M_\varphi$ maps $H^1$ into $H^2$; then using, say, 
La Vall\'ee Poussin kernels, one can obtain the right side.)  

\n We do
not know how to characterize the functional
$\|\varphi\|_{[3]}$ is an analogous fashion. However, the
preceding results show:

\proclaim Theorem 2.8. For any $q\ge 1$, let
$$\alpha(q) = \sup\left\{{\|\varphi\|_{[3]}\over \|\varphi\|_{[2]}}
\right\}$$
where the supremum runs over all non-zero functions $\varphi\colon \ \NN\to
\CC$ with support in $[0,q]$. We have then 
$$a_1\sqrt q \le \alpha(2^q) \le a_2 \sqrt q\leqno
(2.4)$$
 where $a_1,a_2$ are positive
absolute constants.

\pf By the preceding equivalent description of
$\|\varphi\|_{[2]}$,
if $\varphi$ is supported by $[0,2^q]$,  we clearly have
$  \|\varphi\|_2  \le a_3 \|\varphi\|_{[2]}\sqrt q$. Recall
that by (1.8) and (1.2) 
$\|\varphi\|_{M_3(\NN)}\le \|\varphi\|_{B_1(\NN)}  \le
\|\varphi\|_2$. Hence 
 $
\|\varphi\|_{M_3(\NN)}  \le a_3\sqrt q \ \|\varphi\|_{[2]}.  $
This yields the right side of (2.4). The converse follows from
a combination of Lemmas~2.6 and 2.7.\qed

 \bigskip
\n {\bf \S 3.    $M_3(\NN)$ viewed as a space of completely
bounded multipliers}

\n It is known (see [P1, 
p.~109]) that $\|{\varphi}\|_{M_2(\NN)}$ is equal to the $cb$
norm of $M_{\varphi}$ viewed as a  multiplier from $H^1$ to
$H^1$.
  Analogously, we will now show that $\|{\varphi}\|_{M_3(\NN)}$
coincides  with the c.b.\ norm of the ``shifted multiplier''
defined
 by ${\varphi}$ as above but acting 
from $H^1$ to $H^1(\ell_\infty)$.  A somewhat similar statement
can also be proved for general discrete groups
or semi-groups,
but the next statement uses the commutativity of $\NN$.

To state this result, we need some background on
 operator space theory from 
[P6]. Following [P6], we call ``natural'' the operator space (o.s.\ in short) 
structure on $L_1({{\T}})$ for which the o.s.\ dual of $L_1({\T})$ is completely 
isometric to $L_\infty({\T})$. Since $H^1$ is a subspace 
of $L_1({\T})$, this also 
induces a ``natural'' o.s.\ structure on $H^1$. More generally, for any o.s.\ 
$E$, we equip the space $L_1({\T}; E)$ with the o.s.\ structure defined by the 
``o.s.\ projective'' tensor product
 $L_1({\T})\otimes^{\wedge} E$ 
introduced by 
Effros-Ruan [ER1--2] and Blecher-Paulsen [BP1], (see also [P6]). As a Banach 
space $L_1({\T}; E)$ is the same as the classical projective tensor product 
$L_1({\T})\buildrel {\wedge}\over
{\otimes} E$ in Grothendieck's sense, but the o.s.\ structure 
encodes additional information that the norm alone does not carry. We define 
$H^1(E)$ as the closed subspace generated by $H^1\otimes E$ (or 
$\hbox{span}(z^n)\otimes E$) in $L_1({\T}; E) = L_1({\T})
\otimes^{\wedge}
E$. 
Again, we will call ``natural'' the o.s.\ structure induced by the embedding 
$H^1(E) \subset L_1({\T}) \otimes^{\wedge} E$. 

For any bounded function $\varphi\colon \ \NN\to \CC$, let
$$\Phi(n) = \sum_{k\ge 0} \varphi(n+k) e_k \in \ell_\infty$$
where $(e_k)_{k\ge 0}$ denotes the canonical basis of $\ell_\infty$. Then for 
any polynomial $x$ in $H^1$ we have
$$\sup_{k\ge 0} \left|\sum  \hat x(n) z^n\varphi(n+k)\right| =
\left\|\sum \hat x(n) z^n\Phi(n)\right\|_{\ell_\infty}.$$
Moreover, $z\to \sum \hat x(n) z^n\Phi(n)$ can be viewed as an element of 
$H^1\otimes\ell_\infty\subset H^1(\ell_\infty)$. Thus, $\varphi$
defines a mapping
$T_\varphi\colon \ 
\hbox{span}[z^n, n\ge 0]\to H^1(\ell_\infty)$ such that
$$T_\varphi\left(\sum \hat x(n)z^n\right) = 
\sum \hat x(n)z^n  \Phi(n).\leqno(3.1)$$
Note that $\varphi$ is shift-bounded on $H^1$ iff $T_\varphi$ is bounded from 
$H^1$ to $H^1(\ell_\infty)$. As observed in [P1, \S 6], $\varphi\in M_2(\NN)$ 
iff $M_\varphi\colon \ H^1\to H^1$ is completely bounded and
$$\|\varphi\|_{M_2(\NN)} = \|M_\varphi\colon \ H^1\to H^1\|_{cb}.$$
For $M_3(\NN)$, the analogous result is as follows:

\proclaim Theorem 3.1. Let $H^1$ and $H^1(\ell_\infty)$ be equipped with their 
natural o.s.s.\ as above. Let $\varphi\colon \ \NN\to \CC$ be a bounded 
function. Then $\varphi\in M_3(\NN)$ iff $T_\varphi\colon \ H^1\to 
H^1(\ell_\infty)$ is completely bounded. Moreover
$$\|\varphi\|_{M_3(\NN)} = \|T_\varphi\colon \ H^1\to H^1(\ell_\infty)\|_{cb}.$$

\pf Let us write $L_\infty$ instead of $L_\infty({\T})$. Let $S = (H^1)^\bot 
\subset L_\infty$ and let $q\colon \ L_\infty\to L_\infty/S$ be the quotient 
map. As is well known (see [P1, \S 6]) we have a completely isometric embedding 
$j\colon \ L_\infty/S\to B(\ell_2)$ taking $q(\bar z^n)$ to the Hankel operator
$$\gamma_n = \sum_{i+j=n} e_{ij}.$$
Assume $T_\varphi\colon \ H^1\to H^1(\ell_\infty)$ completely bounded.  We have 
clearly a completely contractive (restriction) map $r\colon \ L_\infty 
\otimes_{\rm min} \ell_1 \to (H^1(\ell_\infty))^*$ where here $\ell_1$ is 
equipped with its natural (= maximal) o.s.\ structure. Then let
$$v = j(T_\varphi)^*r\colon \ L_\infty \otimes_{\rm min} \ell_1\to B(\ell_2).$$
We have clearly $\|v\|_{cb} \le \|T_\varphi\|_{cb}$.

\n We can assume that $\ell_1 = \hbox{span}[U_n\mid n\ge 0]$
with
$U_n \in B(H)$, 
$\|U_n\|=1$. By the fundamental factorization of c.b.\ maps (see e.g.\ [Pa1, 
p.~105] or [P1, p.~57]), we can find suitable operators
$\alpha,\beta$ with 
$\|\alpha\| \ \|\beta\| = \|v\|_{cb}$ and a representation 
$$\pi\colon \  L_\infty \otimes_{\rm min} B(H)\to B(\widehat H)$$
such that
$$v(y) = \alpha\pi(y)\beta.\leqno \forall~y\in L_\infty\otimes \ell_1$$
This implies
$$v(\bar z^n \otimes U_k) = \alpha\pi(\bar z^n\otimes U_k)\beta.\leqno 
\forall~n,k\ge 0$$
For any polynomial $x$ in $H^1$, we have
$$\langle r(\bar z^n \otimes U_k), T_\varphi(x)\rangle = \hat x(n)
\langle U_k, 
\Phi(n)\rangle = \hat x(n)\varphi(k+n).$$
Hence we find
$$v(\bar z^n \otimes U_k) = \gamma_n\varphi(k+n).$$
Letting $U = \pi(\bar z\otimes I)$ and $\widehat U_k = \pi(I\otimes U_k)$, we 
obtain 
$$\gamma_n\varphi(k+n) = \alpha U^n\widehat U_k \beta$$
hence denoting by $S$ the classical
unilateral  shift  (note $U$ and
$\widehat U_k$ commute)
$$\varphi( i+k+j) = \langle \alpha U^{i+j} \widehat U_k \beta
e_j, e_i\rangle
 = 
\langle \alpha U^{i} \widehat U_k U^j \beta e_j,  e_i\rangle
=\langle {S^i}^*\alpha U^{i} \widehat U_k U^j \beta S^j e_0,
  e_0\rangle
$$
which clearly shows that $\|\varphi\|_{M_3(\NN)} \le \|\alpha\| \ \|\beta\|$. So 
we conclude
$$\|\varphi\|_{M_3(\NN)} \le \|T_\varphi\|_{cb}.$$
The converse inequality can be proved by an extension of the argument given 
above for Lemma~2.3.  We simply mention that 
the c.b.\ norm of $T_\varphi$ is the norm of the ``same'' multiplier acting from 
$H^1(S_1)$ to $H^1(S_1[\ell_\infty])$. Here $S_1[\ell_\infty] =  
S_1 \otimes^{\wedge}  \ell_\infty$; moreover, if we are given
$g,h$ in the  unit 
ball of $S_2$ (=the Hilbert-Schmidt class), and operators
$y_k$ in the unit ball of
$B(\ell_2)$, then  the sequence of products $(gy_k h)_{k\ge
0}$ defines an element of the unit ball  of $S_1[\ell_\infty]$
and conversely any element of the unit ball 
is of  this form. We leave the details to the reader.
See the proof of Theorem 3.3 below for more on this point.\qed

\n {\bf Remark 3.2.} 
More generally, let $d\ge2$.
We denote by $\otimes_h$ the Haagerup tensor product for
which we refer the reader to  either
[CS],   [BP1] or [P6]. Let
$$L(\infty, d)=(\ell_1\otimes_h\cdots
 \otimes_h \ell_1)^*\quad (d\  {\rm
 times}).$$
Consider $\varphi$ in $M_d(\NN)$. Let 
$\Phi\colon\ \NN \to L(\infty, d-2)$ be defined by
$$\Phi(n)= \sum \varphi(n+i_1+\cdots +i_{d-2}) 
e_{i_1}\otimes \cdots
 \otimes e_{i_{d-2}} \in L(\infty, d-2)$$
(note that this series is meant only in the weak-$*$-sense).
Then,    consider again   the multiplier defined in (3.1)
$$T_\varphi \colon H^1 \to H^1(L(\infty, d-2)).$$
The same argument as above shows that
$$ \| \varphi\| _{M_d(\NN)} = \| T_\varphi\| _{cb} .$$

\medskip

\n {\bf Remarks.} 
\item{(i)} It is easy to check that
$\|\varphi\|_{M_3(\NN)}\le 1$ iff the  mapping (generalized
Schur multiplier) taking $e_{ij}\otimes U_k$ to 
$\varphi(i+k+j) e_{ij}$ extends to a complete contraction from 
$B(\ell_2)\otimes_{\rm min}\ell_1$ to $B(\ell_2)$ (or  from 
$K(\ell_2)\otimes_{\rm min} \ell_1$ to $K(\ell_2)$).
 \item{(ii)} Since $\|T_\varphi\|_{cb} \ge \|T_\varphi\|$,
 Lemma 2.3 is a corollary of  Theorem  3.1.
\item{(iii)} The multivariable Schur products have been studied
in [ER3] mainly  with the group case in  mind. Our arguments can
be viewed as a variation on the same theme.
\item{(iv)} The preceding proof is reminiscent of  the main
point in [J].

\medskip

One unpleasant feature of Theorem~3.1 is that it is not obvious from it that
$M_3(\NN)$ is a Banach algebra for the pointwise product, although this is
clear from the definition of $M_3(\NN)$. However, this ``defect'' is
corrected in the next statement, for which we use the following notation: \
 Let $x = (x_k)_{k\ge 0}$ be a sequence in
$H^1(S_1)$. We write
$|||(x_k)||| <c$ if there are factorizations
$$x_k = gy_kh$$
with $g,h\in H^2(S_2)$
$$y_k\in B(\ell_2)\quad \hbox{and}$$
$$\|g\|_{H^2(S_2)} \sup_k \|y_k\|_{B(\ell_2)}
\|h\|_{H^2(S_2)}<c.$$ Then we set
$$|||(x_k)||| = \inf \{c\mid |||(x_k)|||<c\}.$$
This norm is nothing but an
explicit description of the norm in the space
$H^1\otimes^\wedge S_1\otimes^\wedge\ell_\infty$.

\proclaim Theorem 3.3. For any $\varphi\colon \
\NN\to
\CC$ we define
$$|||\varphi|||_{{\cal M}_3} =
\sup\{|||({\bar z}^k M_\varphi(z^kx_k))_{k\ge 0}|||\}$$ where
the supremum runs over all sequences $(x_k)_{k\ge 0}$ in
$H^1(S_1)$ with
$$|||(x_k)|||<1.$$
We have then
$$\|\varphi\|_{M_3(\NN)} = |||\varphi|||_{{\cal M}_3}.$$

\pf Assume $\|\varphi\|_{M_3(\NN)}<1$ so that (2.3) and (2.2) hold. Let $x_k =
gy_kh$ as above. We have then
$$\eqalign{M_\varphi(z^kx_k) &= \sum_{n\ge 0} z^{k+n} \varphi(k+n)
(gy_kh)^\wedge (n)\cr
&= \sum_{i,j} z^k z^i z^j\langle
\xi_i,T_k\eta_j\rangle g_iy_kh_j\cr &=
\left(\sum_i \xi_i\otimes g_iz^i\right)
(z^k T_k\otimes y_k) \left(\sum_j
\eta_j \otimes z^jh_j\right).}$$
Hence
$$\eqalign{&|||({\bar z}^k M_\varphi(z^kx_k))|||\cr
&\le \left\|\sum \xi_i\otimes g_iz^i\right\|
_{H^2({\ell_2}^* \otimes_2 S_2)}
\sup_k\|T_k\otimes y_k\| \left\|\sum \eta_j\otimes
z^jh_j\right\|_{H^2(\ell_2\otimes_2 S_2)}<1.}$$
This shows that $\|\varphi\|_{M_3(\NN)} \ge \|\varphi\|_{{\cal M}_3}$. For the
converse, note that, if $\|\varphi\|_{{\cal M}_3} \le 1$, then a fortiori we
have clearly for any $x$ in the unit ball of $H^1(S_1)$,
$\|(M_\varphi(z^kx))\|_{L^1\otimes^\wedge S^1\otimes^\wedge \ell_\infty}\le 1$
and therefore (see [P6, Lemma~1.7, p.~23]) $\|T_\varphi\colon \ H^1 \to
H^1(\ell_\infty)\|_{cb}\le 1$. Thus we conclude by Theorem~3.1 that
$\|\varphi\|_{M_3(\NN)} \le \|\varphi\|_{{\cal M}_3}$.\medskip

\n {\bf Remark.} The preceding result shows that if we denote by ${\cal
M}_\varphi\colon \ H^1\otimes^\wedge \ell_\infty \to H^1 \otimes^\wedge
\ell_\infty$ the mapping taking
$$z^n\otimes\alpha\qquad  (n\ge 0, \alpha = (\alpha_k)_{k\ge 0} \in
\ell_\infty)$$ to
$$z^n \otimes (\varphi(n+k)\alpha_k)_{k\ge 0},$$
then we have
$$\|\varphi\|_{M_3(\NN)} = \|{\cal M}_\varphi\|_{cb}.$$

 \n {\bf  Remark.} At this point, we are unable  to
prove  that $M_3(\NN)\not=M_4(\NN)$. The main difficulty is
the lack of a "good'' sufficient condition for
$\varphi\in M_3(\NN)$ analogous to Lemma 2.4.

 \n {\bf Final Remarks on shift-boundedness.}
\def\p{\varphi}
The notion of shift-bounded multiplier
obviously makes sense   also for multipliers
from $H^p$ to $H^q$ ($0<p,q\le\infty$). More
generally, if $M_\p$ is a multiplier which is
bounded from
$L_p(\T)$ to $L_q(\T)$, we will say that
it is shift-bounded if, for any $x$ in $L_p(\T)$,
the "maximal function'' (two sided this time)
$$\sup_{k\in \ZZ} |M_\p (z^k x)| $$
is in $L_q(\T)$. This definition has an obvious
extension to more general function spaces than
$L_p(\T)$ and $L_q(\T)$, for instance it makes
sense also if   $L_q(\T)$ is replaced by
the so-called  weak-$L_q$ space which we
will denote by $L_{q,\infty}(\T)$. 

\n By the Nikishin-Maurey theorems (see [M]), any
shift bounded multiplier from
 $L_1(\T)$ to
$L_{1,\infty}(\T)$ is automatically shift bounded
from $L_2(\T)$ to
$L_{2,\infty}(\T)$, and hence by interpolation
it is shift bounded
from $L_p(\T)$ to $L_p(\T)$ for all $1<p<2$. A
similar result holds for multipliers
on the corresponding Hardy spaces.

\n In passing, it is amusing to observe   that
Carleson's celebrated theorem on the a.s.
convergence of Fourier series in $L_2(\T)$ is
essentially equivalent to the assertion that the
Hilbert transform is shift-bounded from $L_2(\T)$
to
$L_{2,\infty}(\T)$ (and actually it is known
[Hu] to be shift-bounded from $L_p(\T)$ to
$L_p(\T)$ for all $1<p<\infty$).

 \bigskip\bigskip
 
\centerline {\bf References}

 \item{[B]} D. Blecher. A completely 
bounded
characterization of operator algebras.
Math. Ann. 303 (1995) 227-240.

\item{[Bo]}  M.  Bo$\dot{z}$ejko. Littlewood functions, Hankel
multipliers and power
bounded operators on a Hilbert space.	
Colloquium Math.   51 (1987) 35-42.

\item{[BoF]}  M. Bo$\dot{z}$ejko and G.  Fendler. 
Herz-Schur multipliers and completely bounded
multipliers of the Fourier algebra of a locally
compact group.  Boll. Unione Mat. Ital.  (7) 3-A
(1984), 297-302.

 \item{[BP1]} D. Blecher and V. Paulsen. Tensor products of operator spaces.
 J. Funct. Anal. 99 (1991) 262-292.
 
\item{[BP2]} $\underline{\hskip1.5in}$. Explicit constructions 
of universal operator algebras and applictions 
to polynomial factorization. Proc.
Amer. Math. Soc. 112 (1991) 839-850.

\item{[BRS]} D. Blecher, Z. J. Ruan and A. Sinclair.
 A characterization
of operator
algebras. J. Funct. Anal. 89 (1990) 188-201.

\item{[CS]}  E. Christensen and A. Sinclair.   A
survey of completely bounded operators.  Bull.
London Math. Soc.  21 (1989) 417-448.

 \item{[ER1]} E. Effros and Z.J. Ruan. On approximation
properties for operator spaces, International J. Math. 1
(1990) 163-187.

 \item{[ER2]} $\underline{\hskip1.5in}$. A new approach to operator spaces.
 Canadian Math. Bull.
34 (1991) 329-337.

 \item{[ER3]} $\underline{\hskip1.5in}$.
 Multivariable multipliers 
for groups and their
operator algebras.
 Operator theory: operator algebras and
applications, Part 1 (Durham, NH, 1988),
197-218, Proc. Sympos. Pure Math., 51, Part 1,
Amer. Math. Soc., Providence, RI, 1990.

\item{[Ga]} J. Garnett. {\it Bounded analytic
functions.}  Academic Press. New-York 1981.

\item{[H]} U. Haagerup.   $M_{0} A(G)$ functions
which are not coefficients of
uniformly bounded representations.  
    Unpublished manuscript, 1985.

\item{[Hu]} R. A. Hunt. 
On the convergence of Fourier series. 
Orthogonal Expansions and their Continuous
 Analogues (Proc. Conf., Edwardsville, Ill.,
1967) pp. 235--255  Southern Illinois Univ.
Press, Carbondale, Ill, 1968.

\item{[J]}  P. Jolissaint. A characterization of
completely bounded multipliers of Fourier algebras. 
Colloquium Math. 63 (1992) 311-313.

\item{[KLM]} N. Kalton and C. Le Merdy.
Solution of a problem of Peller concerning
similarity. J. Op. Theory. To appear.

\item{[M]} B. Maurey. Nouveaux th\'eor\`emes de Nikishin.
Expos\'es Nos. 4 et
5. S\'eminaire Maurey-Schwartz 1973--1974,
\'Ecole Polytech., Paris, 1974.

\item{[Ni]} N. Nikolskii. {\it Treatise on the shift operator.}
Springer Verlag. Berlin 1986.

  \item{[OP]} T. Oikhberg and G. Pisier.
 The ``maximal" tensor product of
two operator spaces. Proc. Edinburgh Math. Soc. 42 (1999)
267--284.

\item{[Pa1]}   V. Paulsen.  {\it Completely bounded
maps and dilations.}  Pitman Research Notes in
Math. 146, Longman, Wiley, New York, 1986.

\item{[Pe]}   V. Peller. Estimates of functions of
power bounded operators on Hilbert space. J.
Oper. Theory  7 (1982) 341-372.

 \item{[P1]}    G. Pisier.
{\it Similarity problems and completely bounded maps.}
Springer Lecture notes 1618 (1995).

\item{[P2]}    $\underline{\hskip1in}$. Multipliers and lacunary
sets in non amenable groups. {  Amer. J. Math.
117 } (1995) 337-376.

 \item{[P3]} $\underline{\hskip1in}$.  The
similarity degree of an operator algebra. St.
Petersburg Math. J. (1998) To appear.

 \item{[P4]} $\underline{\hskip1in}$.  The
similarity degree of an operator algebra. II.
Math. Zeit. To appear.

 \item{[P5]} $\underline{\hskip1in}$. Are unitarizable groups
amenable? To appear.

\item{[P6]} $\underline{\hskip1in}$. {\it  Noncommutative
 vector valued $L_p$-spaces and completely
$p$-summing maps.}  Soc. Math. France. 
Ast\'erisque 237 (1998). 
 
\item{[Ru]} W. Rudin. {\it Fourier analysis on groups.}
Interscience, Wiley New-York, 1962.

\bye

\bye